\documentclass[11pt, a4paper]{amsart}
\usepackage{amssymb, array, amsmath, amscd, pdfpages, enumerate, amsthm, fixltx2e, setspace, mathtools}
\usepackage[papersize={17cm,26cm},total={12.7cm,22.5cm},left=2.15cm]{geometry}

\usepackage[utf8]{inputenc}
\usepackage{amssymb}
\usepackage{bm}
\usepackage{overpic}
\usepackage{xcolor}
\usepackage{graphicx} 
\usepackage{ifthen}
\usepackage{listings}
\usepackage[colorlinks=true]{hyperref}
\usepackage{hyperref}

\newcommand{\COMMapproxRC}[1]{{\color{gray}#1}}
\renewcommand{\COMMapproxRC}[1]{}

\usepackage{ifthen}
\newcommand*{\keyterm}[1]{\emph{#1}}

\newcommand*{\ldelim}[2]{\csname#1l\endcsname#2}   
\newcommand*{\rdelim}[2]{\csname#1r\endcsname#2}   
\newcommand*{\mdelim}[2]{\csname#1m\endcsname#2}   

\newcommand{\gl}{\lambda}
\newcommand{\gw}{\omega}

\newcommand{\eps}{\varepsilon}

\newcommand*{\C}{{\mathbb{C}}}     

\newcommand*{\R}{{\mathbb{R}}}     
\newcommand*{\N}{{\mathbb{N}}}     

\newcommand*{\norm}[1]{\lVert#1\rVert}

\newcommand*{\set} [1]{\{#1\}}
\newcommand*{\setm}[2]{\{\,#1\mid#2\,\}}   
\newcommand*{\iprod}[2]{\langle#1,#2\rangle}    

\newcommand*{\Abs}[2][default]{\ifthenelse{\equal{#1}{default}}{\left\lvert#2\right\rvert}{\ldelim{#1}{\lvert}#2\rdelim{#1}{\rvert}}}
\newcommand*{\Norm}[2][default]{\ifthenelse{\equal{#1}{default}}{\left\lVert#2\right\rVert}{\ldelim{#1}{\lVert}#2\rdelim{#1}{\rVert}}}

\newcommand*{\Iprod}[3][default]{\ifthenelse{\equal{#1}{default}}{\left\langle#2,#3\right\rangle}{\ldelim{#1}{\langle}#2,#3\rdelim{#1}{\rangle}}}
\newcommand*{\Dualpair}[3][default]{\ifthenelse{\equal{#1}{default}}{\left\langle#2,#3\right\rangle}{\ldelim{#1}{\langle}#2,#3\rdelim{#1}{\rangle}}}

\newcommand*{\List}[2][1]{\set{#1,\ldots,#2}}

\newcommand*{\Set}[2][default]{\ifthenelse{\equal{#1}{default}}{\left\{#2\right\}}{\ldelim{#1}{\{}#2\rdelim{#1}{\}}}}

\newcommand*{\dda}[3][1]{\ifthenelse{\equal{#1}{1}}{\frac{d#3}{d#2}}{\frac{d^{#1}#3}{d#2^{#1}}}}
\newcommand*{\ddb}[2][1]{\ifthenelse{\equal{#1}{1}}{\frac{d}{d#2}}{\frac{d^{#1}}{d#2^{#1}}}}
\newcommand*{\pd}[3][1]{\ifthenelse{\equal{#1}{1}}{\frac{\partial{#2}}{\partial{#3}}}{\frac{\partial^{#1}{#2}}{\partial#3^{#1}}}}
\newcommand*{\pdb}[2][1]{\ifthenelse{\equal{#1}{1}}{\frac{\partial}{\partial{#2}}}{\frac{\partial^{#1}}{\partial#2^{#1}}}}

\newcommand*{\Lp}[1][p]{L^{#1}}

\newcommand{\pmat}[1]{\begin{pmatrix}#1\end{pmatrix}}

\newcommand{\mc}[1]{\mathcal{#1}}

\newcommand{\eq}[1]{\begin{align*}#1\end{align*}}
\newcommand{\eqn}[1]{\begin{align}#1\end{align}}

\newcommand{\hiddenproof}[2]{
\ifthenelse{#1}{Hidden}{#2}
}

\newcommand{\citel}[2]{\cite[#2]{#1}}

\newcommand{\longcomm}[1]{}

\newcommand{\blue}[1]{{\color{blue}#1}}

\usepackage{datetime}
\newdateformat{findate}{\THEDAY.\THEMONTH.\THEYEAR}

\DeclareFixedFont{\ttb}{T1}{txtt}{bx}{n}{9} 
\DeclareFixedFont{\ttm}{T1}{txtt}{m}{n}{9}  
\DeclareFixedFont{\itm}{T1}{txtt}{m}{it}{9}  

\usepackage{color}
\definecolor{deepblue}{rgb}{0,0,0.5}
\definecolor{deepred}{rgb}{0.6,0,0}
\definecolor{deepgreen}{rgb}{0,0.5,0}
\renewcommand{\pmat}[1]{\begin{bmatrix}#1\end{bmatrix}}

\newcommand{\CL}{C_\Lambda}

\theoremstyle{definition}

\numberwithin{equation}{section}

\newenvironment{RORP}{\textbf{The Robust Output Regulation Problem.}\it}{}

\newcommand{\yref}{y_{\mbox{\scriptsize\textit{ref}}}}
\newcommand{\wdist}{w_{\mbox{\scriptsize\textit{dist}}}}

\newcommand{\RORname}{\textbf{RORPack}}
\newcommand{\RORnameplain}{RORPack}

\newcommand{\shorten}[1]{{\color{gray}#1}}
\renewcommand{\shorten}[1]{}

\begin{document}

\title[Introduction to ``\RORnameplain'']{Introduction to ``\RORname''\\[1ex]\footnotesize\mdseries{A Python software package for Robust Output Regulation}}
\thispagestyle{plain}

\author{Lassi Paunonen}
\address{Department of Mathematics, Tampere University, PO.\ Box 692, 33101 Tampere, Finland}
 \email{lassi.paunonen@tuni.fi}

\maketitle

\vspace{-3ex}

\begin{center}
  \today
\end{center}

\vspace{3ex}

\begin{abstract}
  This document contains the mathematical introduction to RORPack --- a Python software library for robust output tracking and disturbance rejection for linear PDE systems. The RORPack library is open-source and freely available at 
https://github.com/lassipau/rorpack/
The package contains functionality for automated construction of robust internal model based controllers, simulation of the controlled systems, visualisation of the results, as well as a collection of example cases on robust output regulation of controlled heat and wave equations.
\end{abstract}

{\small\tableofcontents}

\section{General Description}

\RORname\ is a Python library
for controller design and simulation of robust output tracking and disturbance rejection for linear partial differential equation models. The package contains a number of complete examples on robust controller design and simulation of PDE models
of the following types:
\begin{itemize}
  \item one-dimensional diffusion equations with either boundary or distributed control and observation.
  \item two-dimensional heat equation on a rectangular domain.
  \item one-dimensional wave equations with either distributed or boundary control and observation.
    \COMMapproxRC{
    \item two-dimensional wave equation on an annulus with boundary control and observation
    }
\end{itemize}
New examples will also be added in the future versions of the package.

The purpose of this document is to give a general introduction to the background theory of robust output regulation for linear PDE systems, to describe the mathematical models that are included as the example cases, and to document the usage of the software package on a general level.
 
The purpose of \RORname\ is to serve as a tool to \textit{illustrate} the theory of robust output regulation for distributed parameter systems and it should not (yet) be considered as a serious controller design software. 
The developers of the software do not take any responsibility and are not liable for any damage caused through use of this software.
In addition, at its present stage, the library is not optimized in terms of numerical accuracy. Any helpful comments on potential improvements of the numerical aspects will be greatly appreciated!

The \RORname\ software is distributed as open source under the GNU General Public License version 3 (see \texttt{LICENSE.txt} for detailed license and copyright information) and can freely downloaded at the address
\begin{center}
  \href{https://github.com/lassipau/rorpack/}{https://github.com/lassipau/rorpack/}
\end{center}
  This document is published in \href{https://arxiv.org/}{arXiv.org}.
The website of the project is located at the address (hosted by GitHub Pages)
\begin{center}
  \href{https://lassipau.github.io/rorpack/}{https://lassipau.github.io/rorpack/}
\end{center}
All comments and suggestions for improvements are welcome! These can be sent directly to \texttt{\href{mailto:lassi.paunonen@tuni.fi}{lassi.paunonen@tuni.fi}}.

\subsection{Controller Design and Simulation Workflow}

Basic workflow for robust controller design and simulation for a given system is that the user creates a main file with the following parts:
\begin{enumerate}
  \item Calling of user-defined Python routines that return a numerical approximation of the linear PDE model.
  \item Defining the reference and disturbance signals to be considered.
  \item Calling of a \RORname\ routine for construction of a robust controller that is suitable for the type of PDE model.
    This typically involves choosing appropriate parameters for the controller construction.
  \item Calling of \RORname\ routines for construction and simulation of the closed-loop system. The routines return numerical data describing the output, the error, the control signal, and the closed-loop state.
  \item Calling of \RORname\ routines for visualising the behaviour of the output and the output error, as well as 
    user-defined routines for illustrating the behaviour of the state of the controlled PDE system (multidimensional plots or animations).
\end{enumerate}
The example cases included in \RORname\ are built around main files that follow the same structure, and the users are encouraged to implement their own simulations by modifying the example codes.

It is also possible to combine other Python numerical software packages in the study of robust controller design for PDE models in order to employ ready-made numerical approximations or additional numerical methods in the computation of required controller parameters. This approach is illustrated in one of the PDE example cases ``\texttt{heat\_1d\_2}'' described in Section~\ref{sec:Heat1D2} (Case~2). 
In this example the transfer function of the controlled PDE and other parameters used in the controller construction are computed with the \textbf{Chebfun package}~\cite{DriHal14book,Tre13book} (\href{https://chebfun.org/}{https://chebfun.org/}) which provides flexible tools for solution of differential equations with extreme accuracy using spectral methods. 

\section{Introduction to Robust Output Regulation}

In this section we give a general introduction to the mathematical theory of \keyterm{robust output regulation}.
The purpose of the \RORname\ package is to illustrate controller design for linear distributed parameter systems of the form
\begin{subequations}
  \label{eq:plantintro}
  \eqn{
  \dot{x}(t)&=Ax(t)+Bu(t) + B_d\wdist(t), \qquad x(0)=x_0\in X\\
  y(t)&=\CL x(t)+Du(t) + D_d\wdist(t)
  } 
\end{subequations}
on a Banach or a Hilbert space $X$.
Controlled linear PDE models describing diffusion-convection phenomena, waves and vibrations and elastic deformations can be written in this form with a suitable differential operator~$A$~\cite{CurZwa95book,TucWei09book,JacZwa12book}. 
In our main control problem the goal is to design a dynamic error feedback controller in such a way that the output $y(t)$ of the system converges to a given reference signal $\yref(t)$ despite the external disturbance signal $\wdist(t)$. In addition, the controller needs to be \keyterm{robust} in the sense that it achieves the output tracking and disturbance rejection even if the parameters $(A,B,B_d,C,D,D_d)$ are perturbed or contain small uncertainties.

Our main emphasis is on robust output regulation of \keyterm{diffusion-convection equations}, \keyterm{wave equations}, and \keyterm{beam and plate equations}. However, the \RORname\ package can also be used for construction of controllers for  finite-dimensional systems with given matrices $(A,B,B_d,C,D,D_d)$. 

The full robust output regulation problem is defined in the following way. 

\medskip

\begin{RORP}
  Given a reference signal $\yref(t)$, design a dynamic error feedback controller such that the output $y(t)$ of the system converges to the reference signal asymptotically, i.e.,
  \eqn{
    \label{eq:errconv}
    \lim_{t\to \infty}\; \norm{y(t)-\yref(t)}_Y=0
  }
  despite the disturbance signal $\wdist(t)$.
  Moreover, the controller is required to
be \keyterm{robust} in the sense that it achieves the convergence of the output~\eqref{eq:errconv} even under small uncertainties and changes in the parameters $(A,B,B_d,C,$ $D,D_d)$ of the system.
\end{RORP}

The reference signal $\yref(t)$ and the  disturbance signals $\wdist(t)$ are assumed to be of the form
\begin{subequations}
  \label{eq:yrefwdist}
  \eqn{
    \yref(t) &= a_0^1 + \sum_{k=1}^q (a_k^1 \cos(\gw_k t) + b_k^1 \sin(\gw_k t))\\
    \wdist(t) &= a_0^2 + \sum_{k=1}^q (a_k^2 \cos(\gw_k t) + b_k^2 \sin(\gw_k t))
  }
\end{subequations}
for some \textbf{known frequencies} $\set{\gw_k}_{k=0}^q\subset \R$ with $0=\gw_0<\gw_1<\ldots<\gw_q$ and unknown amplitudes $\set{a_k^j}_{k,j},\set{b_k^j}_{k,j}\subset \R$ (some of which may zero). 

\begin{figure}[h!]
  \begin{center}
    \includegraphics[width=.75\linewidth]{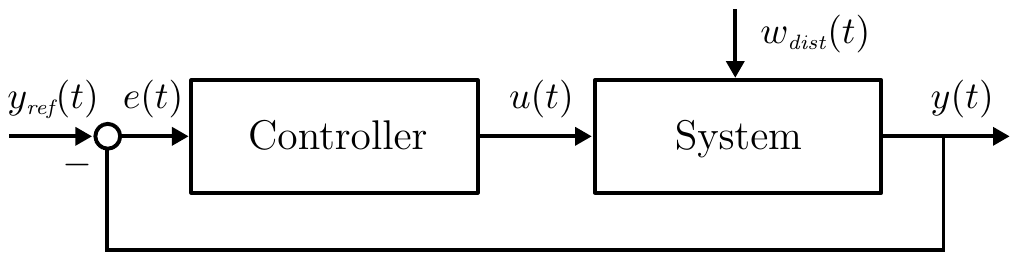} 
    \caption{Dynamic error feedback control scheme.}
    \label{fig:FBcontrol}
  \end{center}
\end{figure}

Dynamic feedback is essential for achieving robust output regulation, and the control problem can indeed be solved with a dynamic error feedback controller (see Figure~\ref{fig:FBcontrol}). The classical \keyterm{internal model principle} gives a characterization for the controllers that solve the robust output regulation problem. This fundamental result was first introduced for finite-dimensional linear systems in the 1970's by Francis and Wonham~\cite{FraWon75a} and Davison~\cite{Dav76} (see~\citel{Hua04book}{Ch. 1} for an excellent overview).
The internal model principle was later extended for infinite-dimensional linear systems by the Systems Theory Research Group at Tampere, Finland in the references~\cite{Imm06phd,PauPoh10,Pau11phd,PauPoh14a,Pau16a}\footnote{Our focus is on linear systems, but there is also an extensive literature on the internal model principle for nonlinear systems, see, e.g.,~\cite{Hua04book,ByrIsi03,ByrIsi04} and references therein.}.

The internal model principle~\cite[Thm. 6.9]{PauPoh10} is also the most important tool in designing controllers for robust output regulation. The result states that a controller solves the robust output regulation problem if and only if the following conditions are satisfied.
\begin{itemize}
  \item The error feedback controller incorporates ``an internal model'' of
    the signals $\yref(t)$ and $\wdist(t)$ in~\eqref{eq:yrefwdist}.
  \item The closed-loop system is exponentially or strongly stable.
\end{itemize}
In controller design, the internal model property can be guaranteed by choosing a suitable \keyterm{structure} for the controller. The rest of parameters are subsequently chosen so that the closed-loop system becomes stable.

The detailed description of the theory of robust output regulation problem and the internal model principle  can be found in the references listed below. The main emphasis in the list is (shamelessly) on the publications by the Systems Theory Research Group at Tampere University, Finland.
\begin{itemize}
  \item \cite{Poh82,XuJer95,HamPoh96a,LogTow97,HamPoh00,RebWei03,LogTow03,ImmPoh06b,Imm07a,DosBas08,BouIdr09,HamPoh10,HamPoh11,PauPoh13a,XuSal14} (including PI-control for PDE models): Robust controller design in various forms for infinite-di\-men\-sional linear systems.
  \item \cite{PauPoh10} and~\cite{PauPoh14a}: The Internal Model Principle for infinite-dimensional linear systems with bounded and unbounded, respectively, input and output operators $B$ and $C$.
  \item \cite{Pau16a,Pau17b}: Robust controller design for \keyterm{regular linear systems}.
  \item \cite{HumPau18,HumKur18}: Robust controller design for \keyterm{port-Hamiltonian systems} and other boundary controlled partial differential equations.
  \item \cite{RebWei03,Pau17carxiv,PauLeGLHMNC18}: Robust controller design \keyterm{impedance passive systems}.
  \item \cite{PauPha18arxiv} for robust finite-dimensional low-order controller design for parabolic systems using Galerkin approximations and model reduction.
  \item References on the output regulation \keyterm{without} the robustness requirement:~\cite{Poh81a,Sch83b,ByrLau00,Deu11,NatGil14,XuDub17a}
\end{itemize}

The considered controllers are of the form
\begin{subequations}
  \label{eq:controller}
  \eqn{
    \dot{z}(t)&=\mc{G}_1  z(t) + \mc{G}_2 (y(t)-\yref(t)),\qquad z(0)=z_0\in Z\\
    u(t) &= K z(t) +D_c (y(t)-\yref(t)).
  }
\end{subequations}
Here $y(t)-\yref(t)$ is the \keyterm{regulation error}.
The construction of the robust controllers are based on the references~\cite{Pau16a,Pau17carxiv} with certain modifications and improvements. In particular, the internal models are defined in their ``real forms'', making the controller real whenever the parameters of the plant are real. The same controllers also achieve robust output tracking also for reference and disturbance signals with complex coefficients $\set{a_k^j}_{k,j},\set{b_k^j}_{k,j}\subset \C$.

The constructions of the robust controllers use the knowledge of the frequencies $\set{\gw_k}_k $ of the reference and disturbance signals, the number of outputs $p:=\dim Y$ of the system~\eqref{eq:plantintro}, and certain knowledge of the system. In particular, the ``minimal controllers'' require knowledge of the values $P(i\gw_k)$ of the transfer function $P(\gl)=\CL R(\gl,A)B+D$ of the system at the complex frequencies $\set{i\gw_k}_k\subset i\R$ of the reference and disturbance signals~\eqref{eq:yrefwdist}.
The other controller structures also use knowledge of the parameters $(A,B,C,D)$ of the system as they involve designing an observer for~\eqref{eq:plantintro}, and require the user to provide stabilizing state feedback and output injection operators. See Section~\ref{sec:ControllerTypes} for more information on the controller-specific requirements.

\section{Using the Software}

The \texttt{rorpack} Python library can be installed for Python 2 or Python 3 as instructed in the \texttt{README.md} file included in the software package (this mainly involves downloading the source codes and typing '\texttt{pip install .}' or '\texttt{pip3 install .}' in the main directory). The subdirectory \texttt{examples/} contains the included PDE examples and simulation files (documented in detail in Section~\ref{sec:PDEcases}).

The following commented example file explains the typical structure and workflow of the controller construction and simulation with \RORname. The considered example cases included in the file \texttt{heat\_1d\_3.py} and documented in Section~\ref{sec:Heat1D3}. Due to the properties of the Python language, the constructor routines used in the main simulation are written at the \textit{beginning}  of the file.

Contents of the file \texttt{heat\_1d\_3.py}: The file begins with comment lines and loading of the necessary parts of the \texttt{rorpack} library as well as other Python packages.

\begin{lstlisting}
'''
Heat equation on the interval [0,1] with  Neumann boundary control
and Dirichlet boundary observation. Approximation with a Finite 
Difference scheme.

Neumann boundary disturbance at x=0, two distributed controls and
two distributed measurements regulated outputs. The controls act 
on the intervals 'IB1' and 'IB2' (Default 'IB1' = [0.3,0.4] and 
'IB2' = [0.6,0.7]) and the measurements are the average 
temperatures on the intervals 'IC1'  and 'IC2' (Default 
'IC1' = [0.1,0.2] and 'IC2' = [0.8,0.9]).
'''
import numpy as np
from rorpack.system import LinearSystem
from rorpack.controller import *
from rorpack.closed_loop_system import ClosedLoopSystem
from rorpack.plotting import *
from laplacian import diffusion_op_1d
\end{lstlisting}

The next part introduces a constructor routine to define the Finite Difference approximation of the heat equation and the input and output operators. The parameter $N$ is the size of the approximation and \texttt{cfun} is a function describing the spatially varying thermal diffusivity of the material.

\begin{lstlisting}
def construct_heat_1d_3(N, cfun, IB1, IB2, IC1, IC2):
    spgrid = np.linspace(0, 1, N+1)

    plt.plot(spgrid,cfun(spgrid))
    plt.title('The thermal diffusivity $c(x)$ of the material')
    plt.tight_layout()
    plt.grid(True)
    plt.show()

    h = spgrid[1]-spgrid[0]
    DiffOp, spgrid = diffusion_op_1d(spgrid, cfun, 'ND')
    A = DiffOp

    B1 = 1/(IB1[1] - IB1[0])*np.logical_and(spgrid >= IB1[0], 
            spgrid <= IB1[1])
    B2 = 1/(IB2[1] - IB2[0])*np.logical_and(spgrid >= IB2[0],
            spgrid <= IB2[1])
    B = np.stack((B1, B2), axis=1)
    C1 = h/(IC1[1] - IC1[0])*np.logical_and(spgrid >= IC1[0],
            spgrid <= IC1[1])
    C2 = h/(IC2[1] - IC2[0])*np.logical_and(spgrid >= IC2[0],
            spgrid <= IC2[1])
    C = np.stack((C1, C2))
    D = np.zeros((2, 2))
    Bd = np.bmat([[np.atleast_2d(-2/h)], [np.zeros((N-1, 1))]])

    return LinearSystem(A, B, C, D, Bd, np.zeros((2,1))), spgrid
 \end{lstlisting}

 The next part defines the parameters of the system and constructs the system $(A,B,B_d,C,D)$ as an object of the class \texttt{LinearSystem} of \RORname.

\begin{lstlisting}
# Parameters for this example.
N = 100

# The spatially varying thermal diffusivity of the material
# cfun = lambda x: np.ones(np.atleast_1d(x).shape)
# cfun = lambda x: 1+x
# cfun = lambda x: 1-2*x*(1-2*x)
cfun = lambda x: 1+.5*np.cos(5/2*np.pi*x)
# Note: Lower diffusivity is difficult for the Low-Gain
# and Passive controllers
# cfun = lambda x: .2-.4*x*(1-x)

# Regions of inputs and outputs
IB1 = np.array([0.3, 0.4])
IB2 = np.array([0.6, 0.7])
IC1 = np.array([0.1, 0.2])
IC2 = np.array([0.8, 0.9])

# Length of the simulation
t_begin = 0
t_end = 8
t_points = 300

# Construct the system.
sys, spgrid = construct_heat_1d_3(N, cfun, IB1, IB2, IC1, IC2)
  \end{lstlisting}

The next part defines the reference signal $\yref(t)$ (in \texttt{yref}) and the disturbance signal $\wdist(t)$ (in \texttt{wdist}) and lists the (real) frequencies $\set{\gw_k}_{k=1}^q$ in the variable \texttt{freqsReal}. Alternative reference and disturbance signals are commented out in the code for further simulation experiments.

\begin{lstlisting}
# Define the reference and disturbance signals, and list the
# required frequencies in 'freqsReal'
# Case 1:
yref = lambda t: np.stack((np.sin(2*t), 2*np.cos(3*t)))
wdist = lambda t: np.sin(6*t)

# Case 2:
# yref = lambda t: np.ones((2,np.atleast_1d(t).shape[0]))
# wdist = lambda t: np.zeros(np.atleast_1d(t).shape)

freqsReal = np.array([0, 1, 2, 3, 6])
\end{lstlisting}

The next part constructs the chosen controller structure and the closed-loop system as objects of the \RORname\ classes. Alternative controller structures are commented out in the code for easy comparison of controller performances.

\begin{lstlisting} 
# Construct the controller and the closed loop system.

# Controller choices, Low-gain robust controller
# Requires the transfer function values P(i*w_k)
# epsgainrange = np.array([0.3,0.6])
# Pvals = np.array(list(map(sys.P, 1j * freqsReal)))
# contr = LowGainRC(sys, freqsReal, epsgainrange, Pvals)


# Dual observer-based controller
# Requires stabilizing operators K and L1
# and the transfer function values P_L(i*w_k)
# K = -sys.B.conj().T
# L1 = -10*sys.C.conj().T
# PLvals = np.array(list(map(lambda freq: sys.P_L(freq, L1),
#     1j * freqsReal)))
# IMstabmargin = 0.5
# IMstabmethod = 'LQR'
# contr = DualObserverBasedRC(sys, freqsReal, PLvals, K, L1, 
#	IMstabmargin, IMstabmethod)

# Observer-based controller
# Requires stabilizing operators K21 and L
# and the transfer function values P_K(i*w_k)
K21 = -sys.B.conj().T
L = -10*sys.C.conj().T
PKvals = np.array(list(map(lambda freq: sys.P_K(freq, K21), 
    1j * freqsReal)))
IMstabmargin = 0.5
IMstabmethod = 'LQR'
contr = ObserverBasedRC(sys, freqsReal, PKvals, K21, L,
	IMstabmargin, IMstabmethod)

# Construct the closed-loop system 
clsys = ClosedLoopSystem(sys, contr) 
\end{lstlisting}

The next part simulates the behaviour of the closed-loop system. The initial state of the system is chosen in the variable \texttt{x0}, and the initial state of the controller is chosen by default to be zero.

\begin{lstlisting}
# Simulate the system, define the initial state x_0
# x0fun = lambda x: np.zeros(np.size(x))
x0fun = lambda x: 0.5 * (1 + np.cos(np.pi * (1 - x)))
# x0fun = lambda x: 3*(1-x)+x
# x0fun = lambda x: 1/2*x**2*(3-2*x)-1
# x0fun = lambda x: 1/2*x**2*(3-2*x)-0.5
# x0fun = lambda x: 1*(1-x)**2*(3-2*(1-x))-1
# x0fun = lambda x: 0.5*(1-x)**2*(3-2*(1-x))-0.5
# x0fun = lambda x: 0.25*(x**3-1.5*x**2)-0.25
# x0fun = lambda x: 0.2*x**2*(3-2*x)-0.5
x0 = x0fun(spgrid)

# z0 is chosen to be zero by default
z0 = np.zeros(contr.G1.shape[0])
xe0 = np.concatenate((x0, z0))

tgrid = np.linspace(t_begin, t_end, t_points)
sol,output,error,control,t = clsys.simulate(xe0,tgrid, yref, wdist)
print('Simulation took %f seconds' % t)
\end{lstlisting}

Finally, the results of the simulation are plotted in separate figures and the behaviour of the controlled state is animated using the user-defined routines.

\begin{lstlisting}
# Plot the output and the error, and animate the behaviour 
# of the controlled state.
plot_output(tgrid, output, yref, 'subplot', 'default')
plot_error_norm(tgrid, error)
plot_control(tgrid, control)

# In plotting and animating the state, fill in the homogeneous
# Dirichlet boundary condition at x=1
sys_state = np.vstack((sol.y[0:N],np.zeros((1,np.size(tgrid)))))
spgrid = np.concatenate((spgrid,np.atleast_1d(1)))

plot_1d_surface(tgrid, spgrid, sys_state, colormap=cm.plasma)
animate_1d_results(spgrid, sys_state, tgrid)
\end{lstlisting}

\section{Implemented Controller Types}
\label{sec:ControllerTypes}

In this section we list the concrete controllers implemented in \RORname.
The documentation of the code includes additional information on the usage of the construction routines.
The controllers are the following:
\begin{itemize}
  \item The ``minimal robust controller'' (including only the internal model), based on references~\cite{HamPoh00},~\citel{Pau16a}{Sec. IV}. Stabilization of the closed-loop system is based on selection of a suitable \keyterm{low-gain parameter} $\eps>0$.\\[1ex]
      Calling sequence for the construction:\\[-1ex]
     \begin{quotation}
       \texttt{LowGainRC(sys,freqsReal,epsgain,Pvals)}
     \end{quotation}
     \medskip
     where \texttt{sys} contains the parameters of the plant, \texttt{freqs} contains the frequencies $\set{\gw_k}_{k=0}^q$ of the signals $\yref(t)$ and $\wdist(t)$ in~\eqref{eq:yrefwdist}, \texttt{epsgain} is the value of the low-gain parameter $\eps>0$. Finally, \texttt{Pvals} is a $(q+1)\times p\times m$ array containing the values $P(i\gw_k)=\CL R(i\gw_k,A)B+D\in \C^{p\times m}$ of the transfer function of~\eqref{eq:plantintro} at the complex frequencies $\set{i\gw_k}_{k=0}^q$ of the reference and disturbance signals in~\eqref{eq:yrefwdist}.
     The parameter \texttt{epsgain} can alternatively be a vector of length 2 providing minimal and maximal values for $\eps$. The controller construction has a naive functionality for finding an $\eps$ to optimize stability margin of the numerically approximated closed-loop system (simply by starting from the minimal value and increasing $\eps$ in steps). 

     \bigskip

  \item The ``observer-based robust controller'', based on references~\citel{HamPoh10}{Sec. 7},~\citel{Pau16a}{Sec. VI},~\citel{Pau17b}{Sec. 5}. The closed-loop stability is achieved using an observer for the state of the system~\eqref{eq:plantintro}.\\[1ex]
      Calling sequence for the construction:\\[-1ex]
     \begin{quotation}
       \texttt{ObserverBasedRC(sys,freqsReal,PKvals,K21,L,\\ 
	 \phantom{a}\hspace{2.7cm} IMstabmargin,IMstabmethod,CKRKvals)}
     \end{quotation}
     \medskip
     where \texttt{sys} contains the parameters of the plant, \texttt{freqsReal} contains the frequencies $\set{\gw_k}_{k=0}^q$ of the signals $\yref(t)$ and $\wdist(t)$ in~\eqref{eq:yrefwdist}.
     The parameters \texttt{K21} and \texttt{L} describe operators $K_{21}$ and $L$, respectively, such that $A+BK_{21}^\Lambda$ and $A+L\CL$ generate exponentially stable semigroups. 
The parameter \texttt{PKvals} is a $(q+1)\times p\times m$ array containing the values $P_K(i\gw_k)=(\CL +DK_{21}^\Lambda)R(i\gw_k,A+BK_{21}^\Lambda)B+D\in \C^{p\times m}$ of the stabilized transfer function of~\eqref{eq:plantintro} at the complex frequencies $\set{i\gw_k}_{k=0}^q$ of the reference and disturbance signals in~\eqref{eq:yrefwdist}.

     Instead of the approach used in~\cite{Pau16a}, the ``internal model'', i.e., the pair $(G_1,B_1)$, in the controller is stabilized using either LQR-based design (\texttt{IMstabmethod = 'LQR'}, by default) or pole placement (\texttt{IMstabmethod = 'poleplacement'}) with a stability margin \texttt{IMstabmargin} (default = \texttt{0.5}). Note that the variable \texttt{IMstabmargin} only determines the stability margin of the internal model, and the stability margin of the closed-loop system also depends on the stability margins of the semigroups generated by $A+BK_2^\Lambda$ and $A+L\CL$.
     Finally, \texttt{CKRKvals} (optional) is a $(q+1)\times p\times N$ array containing elements $(\CL +DK_{21}^\Lambda)R(i\gw_k,A+BK_{21}^\Lambda)\in \C^{p\times N}$ for $k=\List[0]{q}$. If this parameter is not given, the routine computes these values based on the same numerical approximation as the one used in the simulation.
  \item The ``dual observer-based robust controller'', based on references~\citel{Pau16a}{Sec. V},~\citel{Pau17b}{Sec. 4}. The closed-loop stability is achieved using a complementary controller structure that coincides with observer-based stabilization of the dual of the closed-loop system.\\[1ex]
      Calling sequence for the construction:\\[-1ex]
     \begin{quotation}
       \texttt{DualObserverBasedRC(sys,freqsReal,PLvals,K2,L1,\\
	 \phantom{a}\hspace{2.7cm} IMstabmargin,IMstabmethod,RLBLvals)}
     \end{quotation}
     \medskip
     where \texttt{sys} contains the parameters of the plant, \texttt{freqs} contains the frequencies $\set{\gw_k}_{k=0}^q$ of the signals $\yref(t)$ and $\wdist(t)$ in~\eqref{eq:yrefwdist}. The parameters \texttt{K2} and \texttt{L1} describe operators $K_2$ and $L_1$, respectively, such that $A+BK_2^\Lambda$ and $A+L_1\CL$ generate exponentially stable semigroups. 
     The parameter \texttt{PLvals} is a $(q+1)\times p\times m$ array containing the values $P_L(i\gw_k)=\CL R(i\gw_k,A+L_1\CL)(B+L_1D)+D\in \C^{p\times m}$ of the stabilized transfer function of~\eqref{eq:plantintro} at the complex frequencies $\set{i\gw_k}_{k=0}^q$ of the reference and disturbance signals in~\eqref{eq:yrefwdist}.

     Instead of the approach in~\cite{Pau16a}, the ``internal model'', i.e., the pair $(C_1,G_1)$, in the controller is stabilized using either LQR-based design (\texttt{IMstabmethod = 'LQR'}) or pole placement (\texttt{IMstabmethod = 'poleplacement'}) with a predefined stability margin \texttt{IMstabmargin}. Note that the variable \texttt{IMstabmargin} only determines the stability margin of the internal model, and the stability margin of the closed-loop system also depends on the stability margins of the semigroups generated by $A+BK_2^\Lambda$ and $A+L_1\CL$.

     Finally, \texttt{RLBLvals} (optional) is a $(q+1)\times N\times m$ array containing elements $R(i\gw_k,A+L_1\CL)(B+L_1D)\in \C^{N\times m}$ for $k=\List[0]{q}$. If this parameter is not given, the routine computes these values based on the same numerical approximation as the one used in the simulation.

  \item The ``passive minimal controller'' based on~\citel{RebWei03}{Thm. 1.2},~\citel{Pau17carxiv}{Sec. 5.1},~\cite{PauLeGLHMNC18}. For a strongly stable impedance passive system the closed-loop stability can be achieved using passive controller design and a \keyterm{power preserving interconnection} between the passive control system and the controller.\\[1ex]
      Calling sequence for the construction:\\[-1ex]
     \begin{quotation}
	 \texttt{PassiveRC(freqsReal,dim\_Y,epsgain,sys)}
     \end{quotation}
     \medskip
     where
     \texttt{freqsReal} contains the (real) frequencies $\set{\gw_k}_{k=0}^q$ of the signals $\yref(t)$ and $\wdist(t)$ in~\eqref{eq:yrefwdist}, \texttt{dim\_Y} is the dimension of the output space, \texttt{epsgain} is a parameter $\eps>0$ that controls the norm of $K$,
     and \texttt{sys} contains the parameters of the plant.
     In addition, \texttt{Pvals} is a $(q+1)\times p\times m$ array containing the values $P(i\gw_k)=\CL R(i\gw_k,A)B+D\in \C^{p\times m}$ of the transfer function of~\eqref{eq:plantintro} at the complex frequencies $\set{i\gw_k}_{k=0}^q$ of the reference and disturbance signals in~\eqref{eq:yrefwdist}. The parameter \texttt{Pvals} is optional and it is currently not used in the controller design.
     In future versions of the implementation this information may be used to improve the performance of the controller.

The parameter \texttt{epsgain} can alternatively be a vector of length 2 providing minimal and maximal values for $\eps$. The controller construction has a naive functionality for finding an $\eps$ to optimize stability margin of the numerically approximated closed-loop system (simply by starting from the minimal value and increasing $\eps$ in steps). 

It should be noted that the controller construction routine does not test passivity, and for a non-passive system the resulting closed-loop system will typically be unstable.

     \COMMapproxRC{
     \item An ``approximate robust controller'' based on the reference~\cite{HumKur18} for finite-dimensional approximate control design for systems with infinite-dimensional output spaces.\\[1ex]
	 Calling sequence for the construction:\\[-1ex]
	   \begin{quotation}
	     \texttt{ApproximateRC(\dots)}
	   \end{quotation}
	 }
\end{itemize}

\subsection{Comments on the Controller Parameters}

In this section we make some remarks on the choices of the controller parameters.

\textbf{The gain parameter $\eps>0$ in }\texttt{LowGainRC}. 
The theory in~\cite{HamPoh00,RebWei03} guarantees that for a stable system $(A,B,C,D)$ there exists $\eps^\ast>0$ such that for any $0<\eps<\eps^\ast$ the closed-loop system is exponentially stable. The stability margin of the closed-loop system (which directly determines the convergence rate of the regulation error $e(t)$) can be optimized with a suitable choice of $\eps>0$. Finding such an optimal value of $\eps>0$ for a PDE system is a challenging task, but in numerical simulations one can use a naive approach of tracking the spectrum of the finite-dimensional closed-loop system matrix $A_e(\eps)$. This is the approach taken in the example cases in \RORname, though it should be noted that for general PDE systems there is no guarantee that the value $\eps>0$ obtained this way would optimize the closed-loop system for the original PDE system. 
In the example cases we are mainly interested in finding an $\eps>0$ which achieves a reasonable rate of convergence rate of the error $\norm{e(t)}$.

The gain in the passive robust controller \texttt{PassiveRC} can analogously be adjusted with a choice of a parameter $\eps>0$ to optimize the closed-loop stability margin.

\textbf{The transfer function values.} The controllers make use of the values $P(i\gw_k)$, $P_K(i\gw_k)$ and $P_L(i\gw_k)$ of the transfer functions associated to the original or stabilized versions of the PDE system.
The \RORname\ class \texttt{LinearSystem} has the necessary routines for computing these values based on the matrices \texttt{A}, \texttt{B}, \texttt{C}, and \texttt{D}. However, using the same approximation of the PDE for both controller design and simulation corresponds to essentially controlling the approximation as a finite-dimensional system. If possible, to avoid unrealistically positive results, it is therefore better to use two different approximations for controller construction and simulation. This is especially the case if the validity of the approximation for the computation of the parameters can not be guaranteed with absolute certainty. However, there are of course cases, such Galerkin approximations of parabolic systems~\cite{Mor94}, where the values of the controller parameters can be shown to converge with the order of the approximation, and above concerns are unnecessary.
The same comments apply to the parameters \texttt{CKRKvals} and \texttt{RLBLvals}.

For some special PDEs, such as 1D heat or wave equations with constant coefficients, the values $P(i\gw_k)$ may have explicit expressions, but these are very limited special cases. For PDEs with spatially varying parameters there are powerful computational methods that can be used to determine $P(i\gw_k)$, such as the \textbf{Chebfun} package (\href{https://chebfun.org/}{https://chebfun.org/}) employed in \texttt{heat\_1d\_2.py} in Section~\ref{sec:Heat1D2} (``Case 2'').

\section{PDE Models of the Example Cases}
\label{sec:PDEcases}

In this section we introduce the PDE models considered in the example cases and review their fundamental properties.

The considered reference and disturbance signals are in each case combinations of trigonometric functions of the form~\eqref{eq:yrefwdist} with a given set of frequencies. The precise choices of the signals can be seen from the main files of the examples. Similarly, the chosen initial states are visible from the source files, and in several files alternative initial states are provided (these can be used by uncommenting the corresponding lines of code).

\subsection{The 1D Heat Equations}

This collection of examples consider 1D heat equations with spatially varying thermal diffusivity on $\Omega = [0,1]$ with different configurations of control inputs and measured outputs. The first two cases consider control, disturbance, and observation on the boundary, and the final example considers a system with two distributed inputs and outputs. In general, any combination of the above types of inputs and outputs is possible. 
The main property from the point of view of robust output regulation is whether or not the uncontrolled system is exponentially stable (minimal low-gain controller can be used) or impedance passive (the passive robust controller can be used).

In each of the cases, the semidiscretization of the PDE is completed using Finite Differences.

\medskip

\subsubsection*{Case 1: Neumann boundary input at $\xi = 0$, disturbance and output at $\xi = 1$}
~\\[-1ex]

  \noindent Main file name: \texttt{heat\_1d\_1.py}

\medskip

The model on $\Omega = [0,1]$ is
\eq{
  \pd{x}{t}(\xi,t) &=  \pdb{\xi}(c(\xi)\pd{x}{\xi})(\xi,t), \qquad x(\xi,0)=x_0(\xi) \\
      -\pd{x}{\xi}(0,t) &= u(t), \qquad 
      \pd{x}{\xi}(1,t) = \wdist(t), \\
      y(t) &= x(1,t),
    }
    where $c(\cdot)\geq c_0>0$ is the spatially varying \textit{thermal diffusivity} of the material.
    The uncontrolled system is unstable due to the eigenvalue $\gl=0$. 
The example uses the ``Observer-Based Robust Controller'' and ``Dual Observer-Based Robust Controller'' to achieve robust output tracking and disturbance rejection.

\medskip

\subsubsection*{Case 2: Input, output, and disturbance at $\xi = 0$, Dirichlet at $\xi = 1$}
\label{sec:Heat1D2}
~\\[-1ex]

  \noindent Main file name: \texttt{heat\_1d\_2.py}

\medskip

The model on $\Omega = [0,1]$ is
\eq{
  \pd{x}{t}(\xi,t) &=  \pdb{\xi}(c(\xi)\pd{x}{\xi})(\xi,t), \qquad x(\xi,0)=x_0(\xi) \\
      -\pd{x}{\xi}(0,t) &= u(t)+\wdist(t), \qquad 
      x(1,t) = 0, \\
      y(t) &= x(0,t),
    }
    where $c(\cdot)\geq c_0>0$ is the spatially varying thermal diffusivity of the material.
The uncontrolled system is exponentially stable due to the homogeneous Dirichlet boundary condition at $\xi = 1$. The system is also impedance passive since the control input and measured output are collocated,  and because of this the robust output regulation problem can be solved using the Passive Robust Controller.

In this example we use the \textbf{\mbox{Chebfun} package} (\href{https://chebfun.org/}{https://chebfun.org/})~\cite{DriHal14book,Tre13book} in computing the values $P(i\gw_k)$ of the transfer function as well as other parameters required in the controller construction. The Chebfun package utilizes spectral methods and provides powerful and easy-to-use tools for the solution of (especially 1D) boundary value problems with accuracies close to machine precision.
Because of this, Chebfun has great potential in controller design for PDEs and it is especially perfectly suited for robust output regulation of this class of systems.
The only drawback from \RORname's perspective is that at the moment Chebfun is only implemented in Matlab\footnote{Though partial Python implementations exist, the most important functionality for solving BVPs can only be found in the original Matlab version.}. However, using the Python-to-Matlab interface library ``\texttt{matlab.engine}'' included in Matlab (since R2014b), it is possible to call Matlab scripts and functions directly from Python. This approach is used by default in the current PDE example case (this requires installed versions of Matlab, Chebfun package, and the separate installation \texttt{matlab.engine} Python package, see \href{https://se.mathworks.com/help/matlab/matlab-engine-for-python.html}{https://se.mathworks.com/help/matlab/matlab-engine-for-python.html} for details). 
The downside of this approach is that the startup time of the Matlab engine in Python can be extremely slow, and therefore the computations require a relatively long time. However, we get the benefit of very accurate values of the controlled PDE system to be used in the controller design.

Alternate computations of the controller parameters using the Finite Difference approximation (the one used in the main simulation) are commented out in the example code, and can be uncommented to run the \texttt{heat\_1d\_2} example case without the Matlab interface.

\medskip

\subsubsection*{Case 3: Distributed input and output, boundary disturbance at $\xi = 0$, Dirichlet at $\xi = 1$}
\label{sec:Heat1D3}
~\\[-1ex]

  \noindent Main file name: \texttt{heat\_1d\_3.py}

\medskip

The model on $\Omega = [0,1]$ is
\eq{
  \pd{x}{t}(\xi,t) &=  \pdb{\xi}(c(\xi)\pd{x}{\xi})(\xi,t)+b_1(\xi)u_1(t) + b_2(\xi)u_2(t)  \\
      -\pd{x}{\xi}(0,t) &= \wdist(t), \qquad 
      x(1,t) = 0,  \quad x(\xi,0)=x_0(\xi),\\
      y(t) &= \int_0^1\pmat{c_1(\xi)x(\xi,t)\\c_2(\xi)x(\xi,t)}d\xi\in \R^2
    }
    where $c(\cdot)\geq c_0>0$ is the spatially varying thermal diffusivity of the material and
\eq{
  b_1(\xi)&=10\chi_{[.3,.6]}(\xi) , \quad
   b_2(\xi)=10\chi_{[.6,.7]}(\xi), \\
  c_1(\xi)&=10\chi_{[.1,.2]}(\xi), \quad c_2(\xi)=10\chi_{[.8,.9]}(\xi) .
}
    Here $\chi_{[a,b]}(\cdot)$ is the characteristic function on the interval $[a,b]\subset [0,1]$, and thus the control inputs act on the intervals $[0.3,0.4]$ and $[0.6,0.7]$, and the outputs measure the average temperatures on the intervals $[0.1,0.2]$ and $[0.8,0.9]$.
The uncontrolled system is stable because of the homogeneous Dirichlet boundary condition at $\xi = 1$. 

\begin{figure}[h!]
  \begin{minipage}{0.48\linewidth}
    \includegraphics[width=\linewidth]{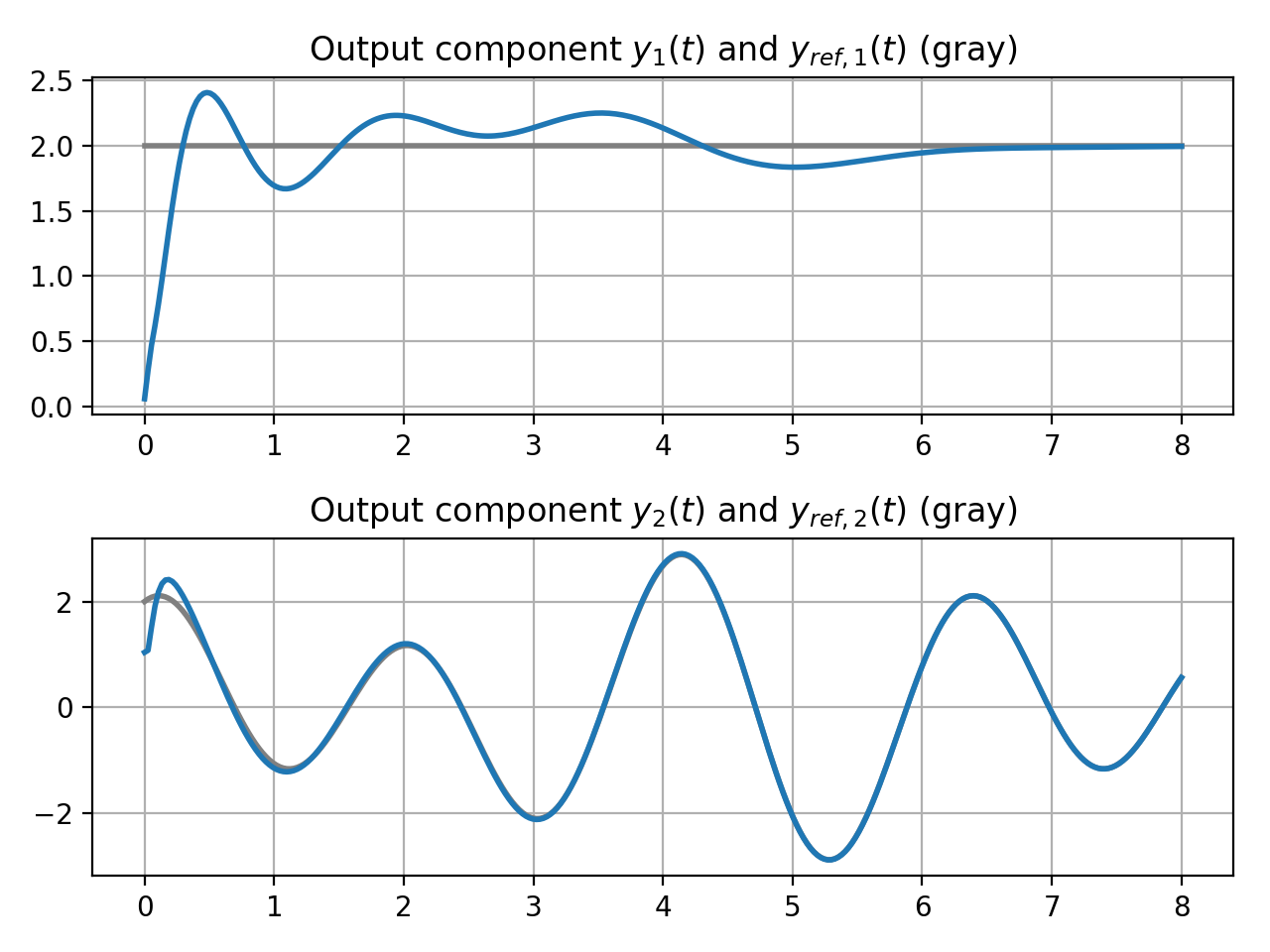}

    \vspace{2ex}

    \includegraphics[width=\linewidth]{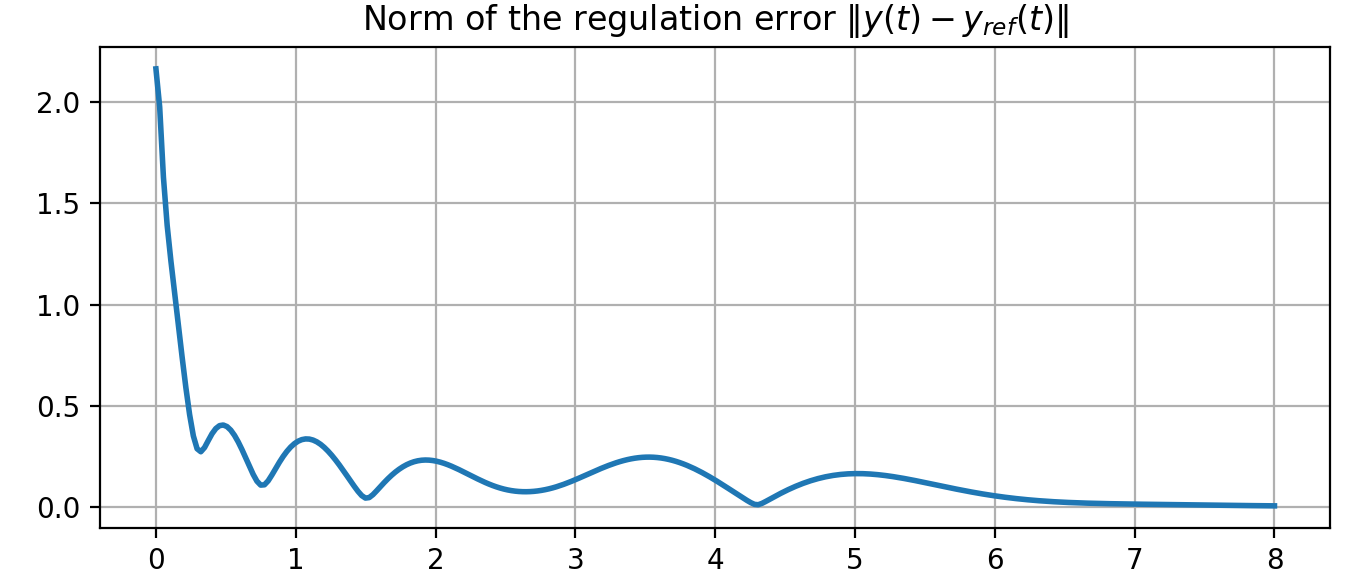}
  \end{minipage}
  \hfill
  \begin{minipage}{0.48\linewidth}
    \includegraphics[width=.95\linewidth]{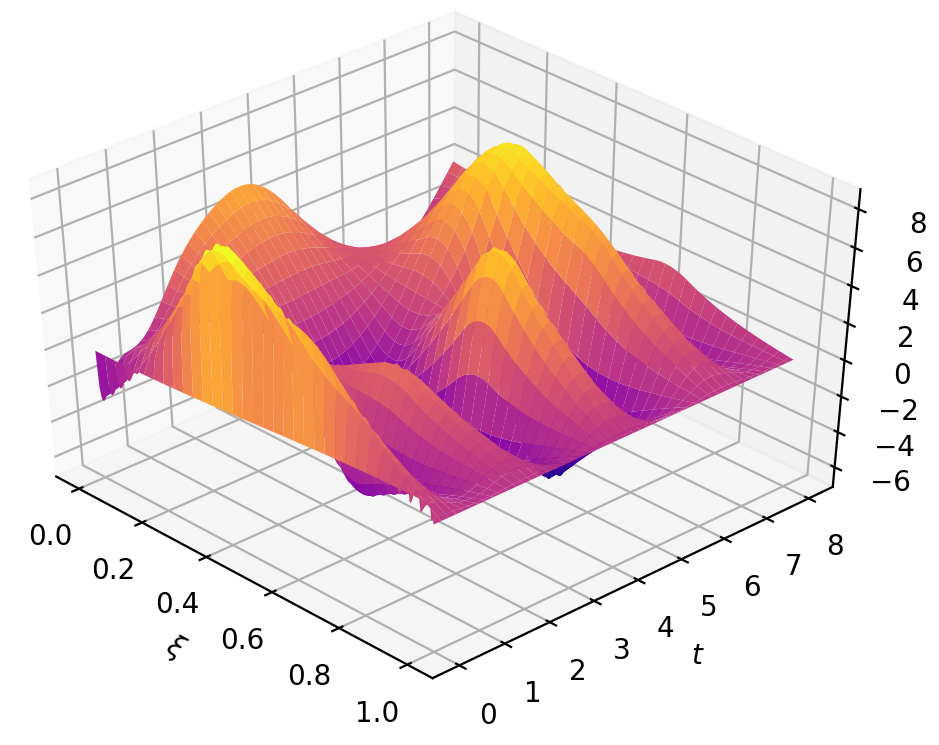}

    \vspace{2ex}

    \begin{center}
      \includegraphics[width=0.93\linewidth]{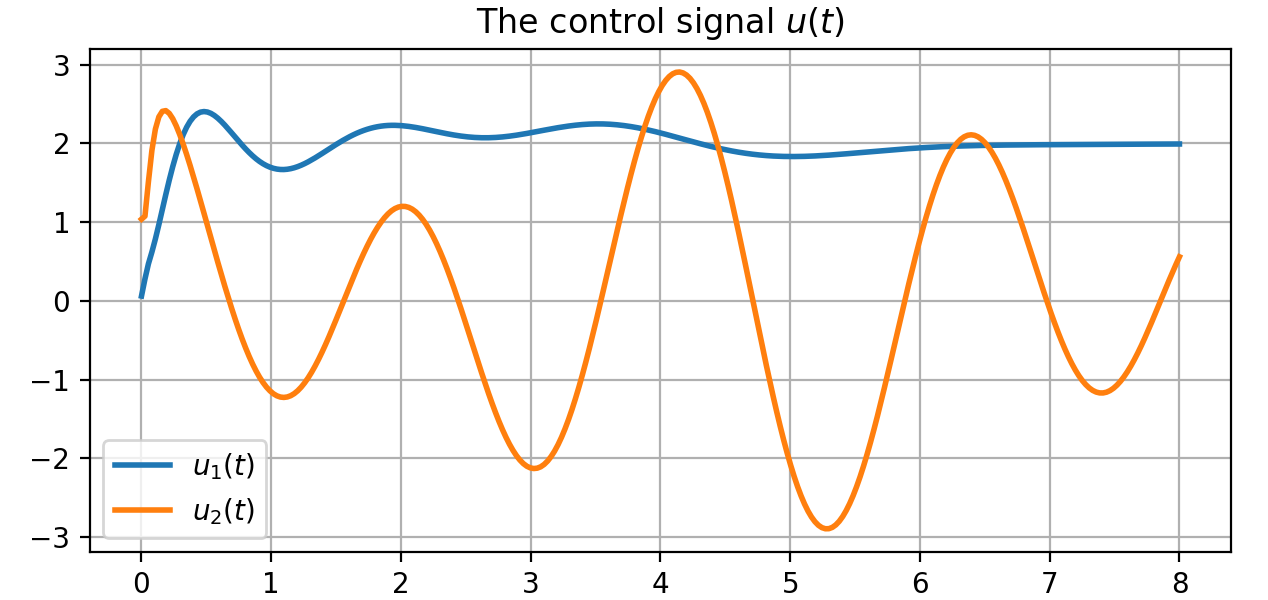}
    \end{center}
  \end{minipage}
\caption{Example output of the 1D Heat equation (``Case 3'').}
  \label{fig:1Dheat3}
\end{figure}

Figure~\ref{fig:1Dheat3} shows example results of the simulations including plots of the outputs and reference signals, norm of the regulation error, computed control signals, and the state of the controlled system as a function of $\xi$ and $t$.

\subsection{The 2D Heat Equations on Rectangular Domains}

These examples consider two-dimensional heat equations on $\Omega = [0,1]\times [0,1]$ with boundary input and output and boundary disturbances. The inputs and outputs act in an averaged sense on the boundaries.

The reference and disturbance signals are again of the form~\eqref{eq:yrefwdist} and can be seen from the main files of the simulations.

\medskip

\subsubsection*{Case 1: Collocated input and output}
~\\[-1ex]

  \noindent Main file name: \texttt{heat\_2d\_1.py}

\medskip

This case is considered in~\citel{Pau16a}{Sec. VII}.
The system on $\Omega = [0,1]\times [0,1]$ with two inputs $u(t)=(u_1(t),u_2(t))^T$ and two outputs $y(t)=(y_1(t),y_2(t))^T$ is determined by
    \eq{
      x_t(\xi,t) &= \Delta x(\xi,t), \qquad x(\xi,0)=x_0(\xi) \\
      \pd{x}{n}(\xi,t)\vert_{\Gamma_1} &= u_1(t), \qquad 
      \pd{x}{n}(\xi,t)\vert_{\Gamma_2} = u_2(t), \qquad
      \pd{x}{n}(\xi,t)\vert_{\Gamma_0} = 0\\
      y_1(t) &= \int_{\Gamma_1}x(\xi,t)d\xi, \qquad
      y_2(t) = \int_{\Gamma_2}x(\xi,t)d\xi.
    }
  Here the parts $\Gamma_0$, $\Gamma_1$, and $\Gamma_2$ of the boundary $\partial \Omega$ are defined so that
  $\Gamma_1 = \setm{\xi=(\xi_1,0)}{0\leq \xi_1\leq 1/2}$,
  $\Gamma_2 = \setm{\xi=(\xi_1,1)}{1/2\leq \xi_1\leq 1}$, 
  $\Gamma_0 = \partial \Omega \setminus (\Gamma_1 \cup \Gamma_2)$. By~\cite[Cor. 2]{ByrGil02} the heat equation defines a regular linear system with feedthrough $D=0$. The system is also impedance passive.
The uncontrolled system is unstable due to the eigenvalue $\gl=0$, but it can be stabilized exponentially with negative output feedback $u(t) = -\kappa y(t)$ for any $\kappa>0$.
In this example we assume the system is pre-stabilized with $\kappa = 1$.

In the simulations, the system is approximated using the eigenmodes of the Laplacian. The system can be controlled either with the low-gain minimal robust controller or the passive controller (with pre-stabilizing negative output feedback), or with one of the observer-based robust controllers.

\medskip

\subsubsection*{Case 2: Non-collocated input and output}
~\\[-1ex]

  \noindent Main file name: \texttt{heat\_2d\_2.py}

\medskip

A similar case (but with collocated inputs and outputs) was considered in~\citel{Pau17carxiv}{Sec. 6.3}.
The system on $\Omega = [0,1]\times [0,1]$ is
    \eq{
      x_t(\xi,t) &= \Delta x(\xi,t), \qquad x(\xi,0)=x_0(\xi) \\
      \pd{x}{n}(\xi,t)\vert_{\Gamma_1} &= u(t), \qquad 
      \pd{x}{n}(\xi,t)\vert_{\Gamma_3} = \wdist(t), \qquad
      \pd{x}{n}(\xi,t)\vert_{\Gamma_0} = 0\\
      y(t) &= \int_{\Gamma_2}x(\xi,t)d\xi.
    }
  Here the parts $\Gamma_0$, $\Gamma_1$, $\Gamma_2$, and $\Gamma_3$ of the boundary $\partial \Omega$ are defined so that
  $\Gamma_1 = \setm{\xi=(0,\xi_2)}{0\leq \xi_2\leq 1}$,
  $\Gamma_2 = \setm{\xi=(\xi_1,1)}{0\leq \xi_1\leq 1}$, 
  $\Gamma_3 = \setm{\xi=(\xi_1,0)}{0\leq \xi_1\leq 1/2}$, 
  $\Gamma_0 = \partial \Omega \setminus (\Gamma_1 \cup \Gamma_2\cup \Gamma_3)$. By~\cite[Cor. 2]{ByrGil02} the heat equation defines a regular linear system with feedthrough $D=0$.
The uncontrolled system is unstable due to the eigenvalue $\gl=0$.  

In the simulations, the system is approximated using Finite Differences with a uniform grid. 

\subsection{The 1D Wave Equations}

In these examples we consider  one-dimen\-sional undamped and damped wave equations with control and observation at the boundaries and inside the domain. 

\medskip

\subsubsection*{Case 1: Non-collocated boundary input and output, disturbance near the output}
~\\[-1ex]

  \noindent Main file name: \texttt{wave\_1d\_1.py}
  
  \medskip

The model on $\Omega = [0,1]$ is
\eq{
  x_{tt}(\xi,t) &=  x_{\xi\xi}(\xi,t), \qquad x(\xi,0)=x_0(\xi), \quad x_t(\xi,0)=x_1(\xi)\\
      -\pd{x}{\xi}(0,t) &= \wdist(t), \qquad 
      \pd{x}{\xi}(1,t) = u(t), \\
      y(t) &= x_t(0,t), \qquad y_m(t) = \int_0^1 x(\xi,t)dt
    }
    The uncontrolled system is unstable due to lack of damping. It also cannot be stabilized with output feedback due to the non-collocated configuration of the inputs and outputs.
    The goal is to achieve tracking of the output $y(t)$. The additional measured output $y_m(t)$ is required to achieve closed-loop stability due to the fact that the system is not exponentially detectable with output $y(t)$ (because of the unobservability of the eigenvalue $\gl=0$).

    In the example case, the second output $y_m(t)$ is used to prestabilize eigenvalue $\gl=0$ of the system with preliminary output feedback of the form $u(t)=-\kappa_m y_m(t)+\tilde{u}(t)$ (where $\tilde{u}(t)$ is the new input of the system). There are also other (and better) ways to handle the situation, and these will be implemented in later versions of this example. 
    The pairs $(A,B)$ and $(C,A)$ are stabilized using collocated designs, i.e., we choose $K=-\kappa B^\ast$ and $L=-\ell C^\ast$ to stabilize the pairs $A+BK$ and $A+LC$.
    The main file contains the pre-stabilization gain $K_m>0$ and the gains $\kappa,\ell>0$ as design parameters.

    The example considers output tracking of the velocity $x_t(0,t)$ to arbitrary 2-periodic reference signals. This is achieved by including  frequencies of the form $k\pi$ for $k\in \List{q}$ in the internal model.  The reference signal is defined by defining its profile over one period in the  variable '\texttt{yref1per}'. Note that it is not necessary to compute the amplitudes of the frequency components of $\yref(t)$ (which would be equivalent to finding the Fourier series  expansion of the reference signal).

\begin{figure}[h!]
  \centering
  \begin{minipage}{0.48\linewidth}
    \begin{flushleft}
      \includegraphics[width=\linewidth]{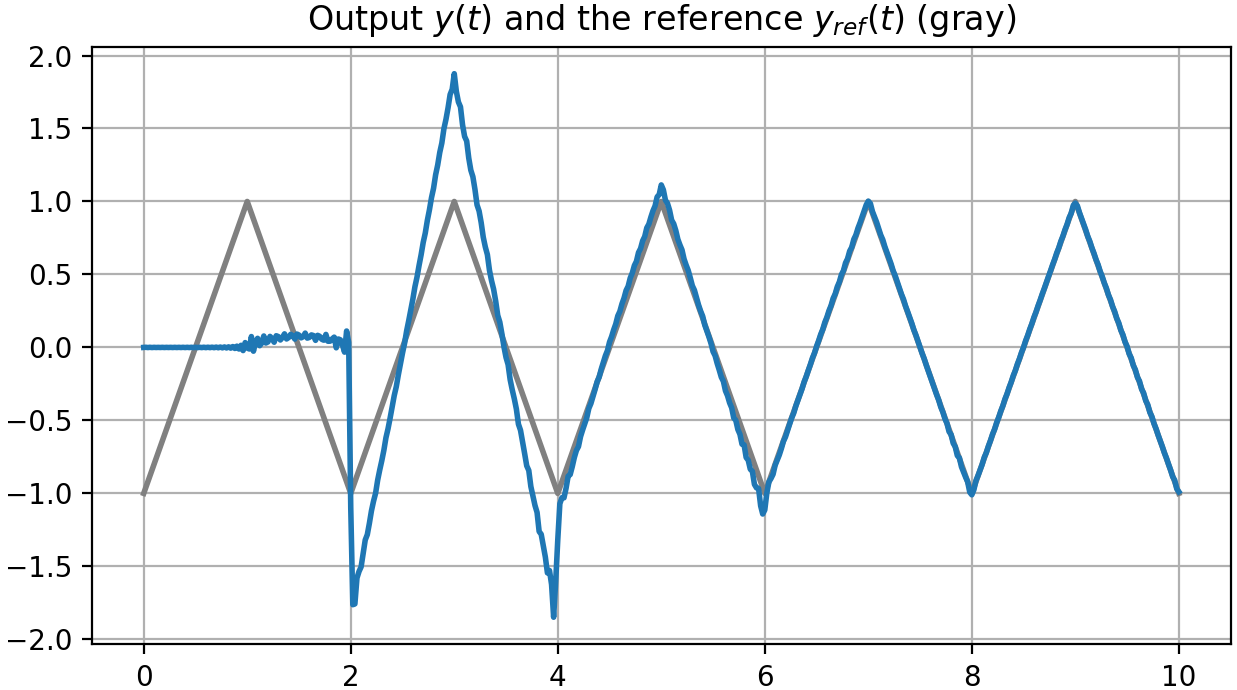}
    \end{flushleft}
  \end{minipage}
  \begin{minipage}{0.48\linewidth}
    \begin{flushright}
      \includegraphics[width=\linewidth]{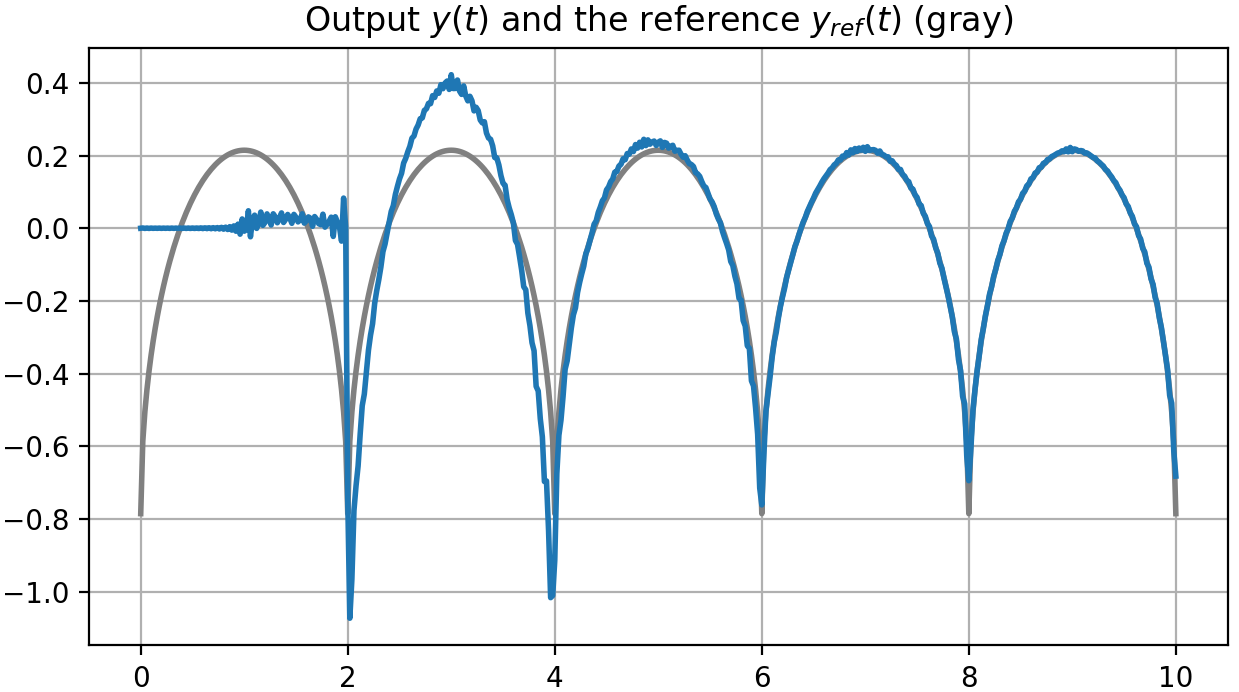}
    \end{flushright}
  \end{minipage}

  \vspace{3ex}

    \begin{center}
      \includegraphics[width=.95\linewidth]{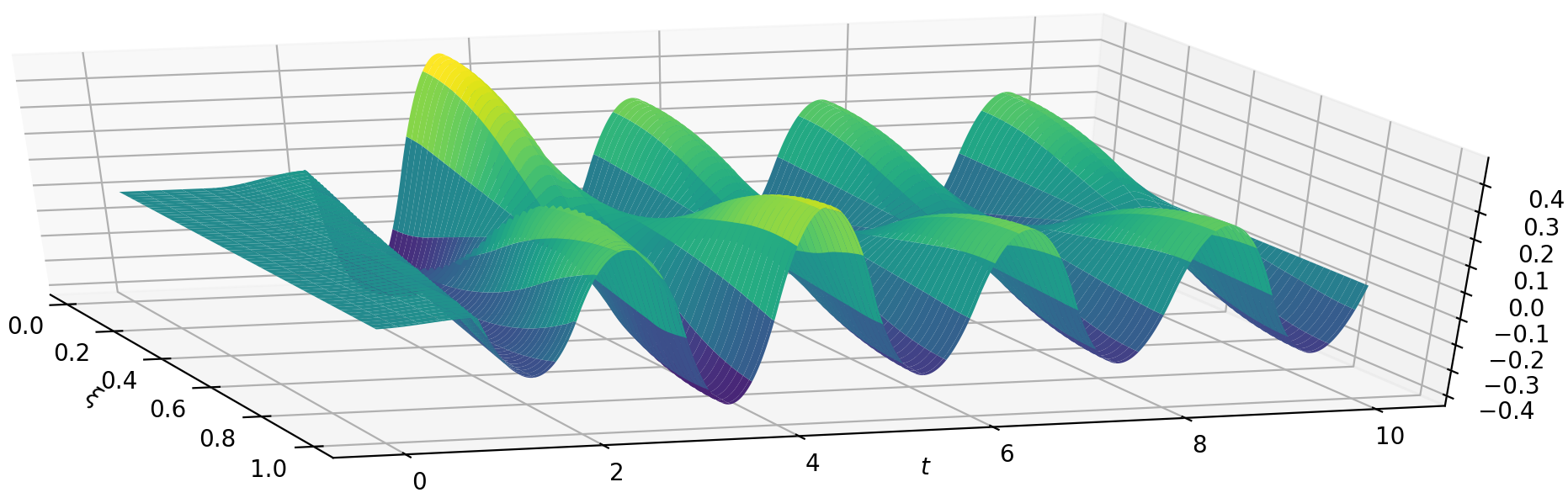}
    \end{center}
    \caption{Example outputs of the periodic tracking for the 1D wave equation (``Case 1'').}
  \label{fig:1Dwave1}
\end{figure}

In the simulations the system is approximated using the orthonormal eigenfunctions of the undamped wave operator.
Figure~\ref{fig:1Dwave1} shows example results of the simulation for two different periodic reference signals --- a nonsmooth triangle signal and signal consisting of semicircles.

\medskip

\subsubsection*{Case 2: Collocated distributed input and output}
~\\[-1ex]

  \noindent Main file name: \texttt{wave\_1d\_2.py}

  \medskip

The model on $\Omega = [0,1]$ is
\eq{
  x_{tt}(\xi,t) &=  x_{\xi\xi}(\xi,t) + b(\xi)u(t) + b_d(\xi)\wdist(t)\\
      x(0,t) &= 0, \qquad 
      x(1,t) = 0 \\
x(\xi,0)&=x_0(\xi), \quad x_t(\xi,0)=x_1(\xi)\\
      y(t) &=  \int_0^1 b(\xi)x(\xi,t)dt
    }
    for some $b(\cdot),b_d(\cdot)\in \Lp[2](0,1;\R)$.
    The uncontrolled system is unstable due to the lack of damping. Since the input and output are distributed, it is also not exponentially stabilizable or detectable, but instead it is only strongly (or polynomially~\cite{BorTom10}) stabilizable provided that $\iprod{b(\cdot)}{\sin(k\pi \cdot)}_{\Lp[2]}\neq 0$ for all $k\in\N$.
    Because of this, it is not possible to use the low-gain controller, and the two observer-based controller designs are not guaranteed to achieve closed-loop stability. However, since the system is impedance passive, the Passive Robust Controller can be used in achieving robust output regulation even in the absence of exponential stability as shown in~\citel{Pau17carxiv}{Sec. 5.1}.
Due to the sub-exponential closed-loop stability, the regulation error does not converge with any uniform convergence rate, but instead the rate depends on the initial state of the system. 

In the simulations the system is approximated using the orthonormal eigenfunctions of the undamped wave operator.

\COMMapproxRC{
  \subsection{The 2D Wave Equation on an Annulus}

  The boundary controlled 2D wave equation we consider is of the form as in~\cite{HumKur18}. In particular, 
  $\Omega\subset\R^n$ is a bounded domain (an open connected set) with a Lipschitz-continuous boundary $\partial\Omega$ split into two parts $\Gamma_0,\Gamma_1$ such that $\overline{\Gamma_0}\cup\overline{\Gamma_1}=\partial\Omega$, $\Gamma_0\cap\Gamma_1=\emptyset$, and $\partial\Gamma_0,\partial\Gamma_1$ both have surface measure zero. The wave equation is of the form
  \eq{
    \rho(\zeta)\frac{\partial^2 w}{\partial t^2} (\zeta,t) &= \nabla\cdot \big(T(\zeta)\,\nabla w\,(\zeta,t)\big), \quad \zeta\in\Omega,\\
    u (\zeta,t) &= \nu\cdot T(\zeta)\,\nabla w(\zeta,t), \quad \zeta\in\Gamma_1,\\
    y (\zeta,t) &= \frac{\partial w}{\partial t}(\zeta,t), \qquad \zeta\in\Gamma_1,\\
    0 &= \frac{\partial w}{\partial t}(\zeta,t),\qquad \zeta\in\Gamma_0,~t>0 \\
    z(\cdot,0) &= w_0,\qquad \frac{\partial z}{\partial t}(\cdot,0)=w_1,
  }
  where $w(\zeta,t)$ is the displacement from the equilibrium at the point $\zeta\in\Omega$ and time $t\geq0$, $\rho(\cdot)$ is the mass density, $T^*(\cdot)=T(\cdot)\in L^2(\Omega;\R^n)$ is Young's modulus and $\nu\in L^\infty(\partial\Omega;\R^n)$ is the unit outward normal at $\partial\Omega$. The functions $\rho(\cdot)$ and $T(\cdot)$ are essentially bounded from both above and below away from zero. In the equation, the boundary input and boundary measurement are collocated.

  In particular, we consider the wave equation on an annulus 
  $\Omega := \set{\zeta \in \mathbb{R}^2 \mid 1 < \|\zeta\| < 2}$. 
  We choose $\partial \Omega = \Gamma_0 \cup \Gamma_1$ where $\Gamma_0 = \set{\zeta \in \partial\Omega \mid \|\zeta\| = 1}$ and $\Gamma_1 = \set{\zeta \in \partial\Omega \mid \|\zeta\| = 2}$ so that the control and the measurement are on the outer boundary of the annulus.

  The output space is infinite-dimensional,  $Y :=  \setm{h\in H^{1/2}(\partial\Omega)}{ h\vert_{\Gamma_0}=0}$ which is continuously and densely embedded into $L^2(\Gamma_1)$~\citel{TucWei09book}{Thm.~13.6.10}. 
  The \keyterm{approximate robust controller} is designed based on finite-dimensional Fourier subspaces
  of $Y$, and the corresponding orthonormal projection $P_N$ onto $Y_N$.

  The initial reference and disturbance signals in the example are given by
  \eq{
    y_{ref}(\theta, t)&  = \frac{1}{2\pi^2}(\pi - \theta)^2\sin(\pi t) +
    \frac{1}{2}\sin\left(\frac{\theta}{2}\right)\cos(2\pi t) \\
    w(\theta, t) & = \cos(\theta)\sin(2\pi t) + \sin(\theta)\sin(\pi t).
  }
  The real frequencies are then $\gw_1=\pi$, $ \gw_2=2\pi$.

  The numerical approximation in the example is done using the spectral method based on the eigenvalues and eigenfunctions of the undamped wave equation on the annulus. The approximate controller includes the stabilizing term $D_c$ which corresponds to pre-stabilizing the system with negative output feedback. 
}

\section{Contributors}

\begin{description}
  \item[Lassi Paunonen] Project leader, developer (2017--)
  \item[Mikko Aarnos] Developer (2018)
  \item[Jukka-Pekka Humaloja] Developer (2018--)
\end{description}

  \appendix

\section{External Links}

\begin{itemize}
  \item \href{https://sysgrouptampere.wordpress.com}{https://sysgrouptampere.wordpress.com} --- Systems Theory Research Group at Tampere University.
  \item  \href{https://github.com/lassipau/rorpack/}{https://github.com/lassipau/rorpack/} --- RORPack at GitHub
\item \href{https://lassipau.github.io/rorpack/}{https://lassipau.github.io/rorpack/}  --- RORPack homepage (hosted by GitHub Pages)
  \item \href{http://mathesaurus.sourceforge.net/matlab-numpy.html}{http://mathesaurus.sourceforge.net/matlab-numpy.html} --- Useful information on differences in the Matlab and Python/Numpy syntax for vectors and matrices.
  \item \href{https://chebfun.org}{https://chebfun.org} --- Homepage of the Chebfun Project.
\end{itemize}

\shorten{

  \section{Code: TODO}

  \begin{itemize}
    \item
  \end{itemize}

  TODO [\blue{Later}]
  \begin{itemize}
    \item Integrate the use of \texttt{FEniCS/Dolfin}
    \item Extend construction to systems with complex parameters (need to include values $P(-i\gw_k)$ etc, find a general way to handle both real and complex systems at the same time?
    \item improved color choices for the plots, define ``styles'', e.g., ``temperature'', ``general'', ``deflection''
    \item Youtube Video lecture on usage of RORPack!
  \end{itemize}
}


\begin{thebibliography}{10}

\bibitem{BorTom10}
Alexander Borichev and Yuri Tomilov.
\newblock Optimal polynomial decay of functions and operator semigroups.
\newblock {\em Math. Ann.}, 347(2):455--478, 2010.

\bibitem{BouIdr09}
S.~Boulite, A.~Idrissi, and A.~Ould~Maaloum.
\newblock Robust multivariable {PI}-controllers for linear systems in {B}anach
  state spaces.
\newblock {\em J. Math. Anal. Appl.}, 349(1):90--99, 2009.

\bibitem{ByrGil02}
Christopher~I. Byrnes, David~S. Gilliam, Victor~I. Shubov, and George Weiss.
\newblock Regular linear systems governed by a boundary controlled heat
  equation.
\newblock {\em Journal of Dynamical and Control Systems}, 8(3):341--370, 2002.

\bibitem{ByrIsi03}
Christopher~I. Byrnes and Alberto Isidori.
\newblock Limit sets, zero dynamics, and internal models in the problem of
  nonlinear output regulation.
\newblock {\em IEEE Trans. Automat. Control}, 48(10):1712--1723, 2003.

\bibitem{ByrIsi04}
Christopher~I. Byrnes and Alberto Isidori.
\newblock Nonlinear internal models for output regulation.
\newblock {\em IEEE Trans. Automat. Control}, 49(12):2244--2247, 2004.

\bibitem{ByrLau00}
Christopher~I. Byrnes, Istv{\'a}n~G. Lauk{\'o}, David~S. Gilliam, and Victor~I.
  Shubov.
\newblock Output regulation problem for linear distributed parameter systems.
\newblock {\em IEEE Trans. Automat. Control}, 45(12):2236--2252, 2000.

\bibitem{CurZwa95book}
Ruth~F. Curtain and Hans~J. Zwart.
\newblock {\em An Introduction to Infinite-Dimensional Linear Systems Theory}.
\newblock Springer-Verlag, New York, 1995.

\bibitem{Dav76}
Edward~J. Davison.
\newblock The robust control of a servomechanism problem for linear
  time-invariant multivariable systems.
\newblock {\em IEEE Trans. Automat. Control}, 21(1):25--34, 1976.

\bibitem{Deu11}
Joachim Deutscher.
\newblock Output regulation for linear distributed-parameter systems using
  finite-dimensional dual observers.
\newblock {\em Automatica J. IFAC}, 47(11):2468--2473, 2011.

\bibitem{DosBas08}
V.~Dos~Santos, G.~Bastin, J.-M. Coron, and B.~d'Andr\'{e}a Novel.
\newblock Boundary control with integral action for hyperbolic systems of
  conservation laws: stability and experiments.
\newblock {\em Automatica J. IFAC}, 44(5):1310--1318, 2008.

\bibitem{DriHal14book}
Tobin~A. Driscoll, Nicholas Hale, and Lloyd~N. Trefethen.
\newblock {\em Chebfun Guide}.
\newblock Pafnuty Publications, 2014.

\bibitem{FraWon75a}
Bruce~A. Francis and W.~Murray Wonham.
\newblock The internal model principle for linear multivariable regulators.
\newblock {\em Appl. Math. Optim.}, 2(2):170--194, 1975.

\bibitem{HamPoh96a}
Timo H{\"a}m{\"a}l{\"a}inen and Seppo Pohjolainen.
\newblock Robust control and tuning problem for distributed parameter systems.
\newblock {\em Internat. J. Robust Nonlinear Control}, 6(5):479--500, 1996.

\bibitem{HamPoh00}
Timo H{\"a}m{\"a}l{\"a}inen and Seppo Pohjolainen.
\newblock A finite-dimensional robust controller for systems in the
  {CD}-algebra.
\newblock {\em IEEE Trans. Automat. Control}, 45(3):421--431, 2000.

\bibitem{HamPoh10}
Timo H{\"a}m{\"a}l{\"a}inen and Seppo Pohjolainen.
\newblock Robust regulation of distributed parameter systems with
  infinite-dimensional exosystems.
\newblock {\em SIAM J. Control Optim.}, 48(8):4846--4873, 2010.

\bibitem{HamPoh11}
Timo H{{\"a}}m{{\"a}}l{{\"a}}inen and Seppo Pohjolainen.
\newblock A self-tuning robust regulator for infinite-dimensional systems.
\newblock {\em IEEE Trans. Automat. Control}, 56(9):2116--2127, 2011.

\bibitem{Hua04book}
Jie Huang.
\newblock {\em Nonlinear Output Regulation, Theory and Applications}.
\newblock SIAM, Philadelphia, 2004.

\bibitem{HumKur18}
Jukka-Pekka Humaloja, Mikael Kurula, and Lassi Paunonen.
\newblock Robust output regulation for wave equations with boundary control.
\newblock {\em in review}.

\bibitem{HumPau18}
Jukka-Pekka Humaloja and Lassi Paunonen.
\newblock Robust regulation of infinite-dimensional port-{H}amiltonian systems.
\newblock {\em IEEE Transactions on Automatic Control}, 63(5):1480--1486, 2018.

\bibitem{Imm06phd}
Eero Immonen.
\newblock {\em State Space Output Regulation Theory for Infinite-Dimensional
  Linear Systems and Bounded Uniformly Continuous Exogenous Signals}.
\newblock PhD thesis, Tampere University of Technology, 2006.

\bibitem{Imm07a}
Eero Immonen.
\newblock On the internal model structure for infinite-dimensional systems:
  {T}wo common controller types and repetitive control.
\newblock {\em SIAM J. Control Optim.}, 45(6):2065--2093, 2007.

\bibitem{ImmPoh06b}
Eero Immonen and Seppo Pohjolainen.
\newblock Feedback and feedforward output regulation of bounded uniformly
  continuous signals for infinite-dimensional systems.
\newblock {\em SIAM J. Control Optim.}, 45(5):1714--1735, 2006.

\bibitem{JacZwa12book}
Birgit Jacob and Hans Zwart.
\newblock {\em Linear Port-Hamiltonian Systems on Infinite-Dimensional Spaces},
  volume 223 of {\em Operator Theory: Advances and Applications}.
\newblock Birkh\"{a}user, Basel, 2012.

\bibitem{LogTow97}
Hartmut Logemann and Stuart Townley.
\newblock Low-gain control of uncertain regular linear systems.
\newblock {\em SIAM J. Control Optim.}, 35(1):78--116, 1997.

\bibitem{LogTow03}
Hartmut Logemann and Stuart Townley.
\newblock Adaptive low-gain integral control of multivariable well-posed linear
  systems.
\newblock {\em SIAM J. Control Optim.}, 41(6):1722--1732, 2003.

\bibitem{Mor94}
Kirsten Morris.
\newblock Design of finite-dimensional controllers for infinite-dimensional
  systems by approximation.
\newblock {\em J. Math. Systems Estim. Control}, 4(2):30, 1994.

\bibitem{NatGil14}
V.~Natarajan, D.S. Gilliam, and G.~Weiss.
\newblock The state feedback regulator problem for regular linear systems.
\newblock {\em Automatic Control, IEEE Transactions on}, 59(10):2708--2723,
  2014.

\bibitem{PauPha18arxiv}
L.~{Paunonen}, D.~{Phan}, and P.~{Laakkonen}.
\newblock {Reduced Order Controller Design for Robust Output Regulation of
  Parabolic Systems}.
\newblock {\em ArXiv e-prints}, November 2018.

\bibitem{Pau11phd}
Lassi Paunonen.
\newblock {\em Output Regulation Theory for Linear Systems with
  Infinite-Dimensional and Periodic Exosystems}.
\newblock PhD thesis, Tampere University of Technology, 2011.

\bibitem{Pau16a}
Lassi Paunonen.
\newblock Controller design for robust output regulation of regular linear
  systems.
\newblock {\em IEEE Trans. Automat. Control}, 61(10):2974--2986, 2016.

\bibitem{Pau17b}
Lassi Paunonen.
\newblock Robust controllers for regular linear systems with
  infinite-dimensional exosystems.
\newblock {\em SIAM J. Control Optim.}, 55(3):1567--1597, 2017.

\bibitem{Pau17carxiv}
Lassi Paunonen.
\newblock {Stability and Robust Regulation of Passive Linear Systems}.
\newblock {\em ArXiv e-prints (https://arxiv.org/abs/1706.03224)}, June 2017.

\bibitem{PauLeGLHMNC18}
Lassi Paunonen, Yann~Le Gorrec, and H\'{e}ctor Ram{\'{\i}}rez.
\newblock A simple robust controller for port--{H}amiltonian systems.
\newblock In {\em Proceedings of the 6th IFAC Workshop on Lagrangian and
  Hamiltonian Methods for Nonlinear Control}, Valparaiso, Chile, May 1--4 2018.

\bibitem{PauPoh10}
Lassi Paunonen and Seppo Pohjolainen.
\newblock Internal model theory for distributed parameter systems.
\newblock {\em SIAM J. Control Optim.}, 48(7):4753--4775, 2010.

\bibitem{PauPoh13a}
Lassi Paunonen and Seppo Pohjolainen.
\newblock Reduced order internal models in robust output regulation.
\newblock {\em IEEE Trans. Automat. Control}, 58(9):2307--2318, 2013.

\bibitem{PauPoh14a}
Lassi Paunonen and Seppo Pohjolainen.
\newblock The internal model principle for systems with unbounded control and
  observation.
\newblock {\em SIAM J. Control Optim.}, 52(6):3967--4000, 2014.

\bibitem{Poh81a}
Seppo Pohjolainen.
\newblock A feedforward controller for distributed parameter systems.
\newblock {\em Internat. J. Control}, 34(1):173--184, 1981.

\bibitem{Poh82}
Seppo~A. Pohjolainen.
\newblock Robust multivariable {PI}-controller for infinite-dimensional
  systems.
\newblock {\em IEEE Trans. Automat. Control}, 27(1):17--31, 1982.

\bibitem{RebWei03}
Richard Rebarber and George Weiss.
\newblock Internal model based tracking and disturbance rejection for stable
  well-posed systems.
\newblock {\em Automatica J. IFAC}, 39(9):1555--1569, 2003.

\bibitem{Sch83b}
J.~M. Schumacher.
\newblock Finite-dimensional regulators for a class of infinite-dimensional
  systems.
\newblock {\em Systems Control Lett.}, 3:7--12, 1983.

\bibitem{Tre13book}
Lloyd~N. Trefethen.
\newblock {\em {Approximation Theory and Approximation Practice}}.
\newblock Philadelphia, PA: Society for Industrial and Applied Mathematics
  (SIAM), 2013.

\bibitem{TucWei09book}
M.~Tucsnak and G.~Weiss.
\newblock {\em Observation and Control for Operator Semigroups}.
\newblock Birkh\"auser Basel, 2009.

\bibitem{XuJer95}
Cheng-Zhong Xu and Hamadi Jerbi.
\newblock A robust {PI}-controller for infinite-dimensional systems.
\newblock {\em Internat. J. Control}, 61(1):33--45, 1995.

\bibitem{XuSal14}
Cheng-Zhong Xu and Gauthier Sallet.
\newblock Multivariable boundary {PI} control and regulation of a fluid flow
  system.
\newblock {\em Math. Control Relat. Fields}, 4(4):501--520, 2014.

\bibitem{XuDub17a}
Xiaodong Xu and Stevan Dubljevic.
\newblock Output and error feedback regulator designs for linear
  infinite-dimensional systems.
\newblock {\em Automatica J. IFAC}, 83:170--178, 2017.

\end{thebibliography}
\end{document}